\documentclass[12pt]{article}

\usepackage{amssymb}
\usepackage{amsfonts}
\usepackage{amsmath}
\usepackage{amsthm}
\usepackage[margin=0.80in]{geometry}
\usepackage{enumerate}
\usepackage{amstext}
\usepackage{layout}

\numberwithin{equation}{section}

\newtheorem{theorem}{Theorem}[section]

\newtheorem{lemma}{Lemma}[section]

\newtheorem{remark}{Remark}[section]
\newtheorem{example}{Example}[section]

\begin{document}

\begin{center}
{\bf\large{A Simple Point Estimator of the Power of Moments}}

\vskip 0.3cm

\vskip 0.3cm
\noindent { {\sc Shuhua Chang}\footnote{Shuhua Chang, Research Center for Mathematics and Economics,
Tianjin University of Finance and Economics, Tianjin 300222, China~~e-mail: szhang@tjufe.edu.cn},
{\sc Deli Li}\footnote{Deli Li, Department of Mathematical
Sciences, Lakehead University, Thunder Bay, ON P7B 5E1, Canada~~e-mail: dli@lakeheadu.ca}, 
{\sc Yongcheng Qi}\footnote{Yongcheng Qi, Department of Mathematics and Statistics, University 
of Minnesota Duluth, Duluth, MN 55812, USA~~e-mail: yqi@d.umn.edu}, and 
{\sc Andrew Rosalsky}\footnote{Andrew Rosalsky, Department of Statistics, University of Florida, 
Gainesville, Florida 32611, USA~~e-mail: rosalsky@stat.ufl.edu}\footnote{Corresponding author.}}
\end{center}

\vskip 0.3cm

\begin{abstract}
Let $X$ be an observable random variable with unknown distribution function $F(x) = \mathbb{P}(X \leq x),
- \infty < x < \infty$, and let
\[\
\theta = \sup\left \{ r \geq 0:~ \mathbb{E}|X|^{r} < \infty \right \}.
\]
We call $\theta$ the power of moments of the random variable $X$. Let $X_{1}, X_{2}, ..., X_{n}$ be a random sample 
of size $n$ drawn from $F(\cdot)$. In this paper we propose the following simple point estimator of $\theta$ and investigate 
its asymptotic properties:
\[
\hat{\theta}_{n} = \frac{\log n}{\log \max_{1 \leq k \leq n} |X_{k}|},
\]
where $\log x = \ln(e \vee x), ~- \infty < x < \infty$. In particular, we show that
\[
\hat{\theta}_{n} \rightarrow_{\mathbb{P}} \theta~~\mbox{if and only if}~~ 
\lim_{x \rightarrow \infty} x^{r} \mathbb{P}(|X| > x) = \infty ~~\forall~r > \theta.
\]
This means that, under very reasonable conditions on $F(\cdot)$, $\hat{\theta}_{n}$ 
is actually a consistent estimator of $\theta$. Hypothesis testing for the power of moments
is conducted and, as an application of our main results, the formula for finding the p-value of 
the test is given. In addition, a theoretical application of our main results is provided together with
three illustrative examples.
\end{abstract}

\vskip 0.5cm


\vskip 0.2cm


\section{Motivation}

The motivation of the current work arises from the following problem concerning parameter estimation.
Let $X$ be an observable random variable with unknown distribution function $F(x) = \mathbb{P}(X \leq x),
- \infty < x < \infty$, and let 
\[\
\theta = \sup\left \{ r \geq 0:~ \mathbb{E}|X|^{r} < \infty \right \}.
\]
We call $\theta$ the {\it power of moments} of the random variable $X$. Clearly $\theta$ is a parameter of the distribution 
of the random variable $X$. Now let $X_{1}, X_{2}, ..., X_{n}$ be a random sample of size $n$ drawn from random 
variable $X$; i.e., $X_{1}, X_{2}, ..., X_{n}$ are independent and identically distributed (i.i.d.) random variables 
whose common distribution function is $F(\cdot)$. It is natural to pose the following question: Can we estimate the the 
parameter $\theta$ based on the random sample $X_{1}$, ..., $X_{n}$? 

This is a serious and important problem. For example, if $\theta > 2$ and if the distribution of $X$ is 
nondegerate, then it is clear that $0 < \mbox{Var} X < \infty$ and so by the classical L\'{e}vy 
central limit theorem, the distribution of
\[
\frac{S_{n} - n \mu}{\sqrt{n}}
\]
is approximately normal (for all sufficiently large $n$) with mean $0$ and variance 
$\sigma^{2} = \mbox{Var}X = \mathbb{E}(X - \mu)^{2}$ where $\mu = \mathbb{E}X$. Thus the problem that we are facing is 
how can we conclude with a high degree of confidence that $\theta > 2$.

In this paper we propose the following point estimator of $\theta$ and will investigate its asymptotic properties:
\[
\hat{\theta}_{n} = \frac{\log n}{\log \max_{1 \leq k \leq n} |X_{k}|}.
\]
Here and below $\log x = \ln(e \vee x), ~- \infty < x < \infty$. 

Our main results will be stated in Section 2 and they all pertain to a sequence of i.i.d. random variables 
$\{X_{n}; n \geq 1\}$ drawn from the distribution function $F(\cdot)$ of the random variable $X$. 
The procedure of our study is as follows.

{\bf Step 1}.~~For deducing the asymptotic properties of $\hat{\theta}_{n}$, $n \geq 1$, we will first precisely 
determine the values of $\rho_{1}, \rho_{2} \in [0, \infty]$ such that
\[
\limsup_{n \rightarrow \infty} \frac{\log \max_{1 \leq k \leq n} X_{k}}{\log n} 
= \frac{1}{\rho_{1}}~~\mbox{almost surely (a.s.)~and}~~
\liminf_{n \rightarrow \infty} \frac{\log \max_{1 \leq k \leq n} X_{k}}{\log n} = \frac{1}{\rho_{2}}~~\mbox{a.s.,}
\]
where $1/0 = \infty$ and $1/\infty = 0$; see Theorems 2.1 and 2.2. 

{\bf Step 2}.~~Following from Theorems 2.1 and 2.2, in Theorem 2.3 we will provide different necessary and sufficient conditions for 
\begin{equation}
\lim_{n \rightarrow \infty} \frac{\log \max_{1 \leq k \leq n} X_{k}}{\log n} = \frac{1}{\rho} 
~~\mbox{a.s. for some}~~ \rho \in [0, \infty].
\end{equation}

{\bf Step 3}.~~Under the assumption that (1.1) holds for some $0 < \rho < \infty$, in Theorem 2.4 we will establish large deviation 
probabilities for 
\[
\mathbb{P}\left(\frac{\log \max_{1 \leq k \leq n} X_{k}}{\log n} \geq \frac{1}{\rho} + x \right)~~\forall~ x > 0 ~~\mbox{and}~~
\mathbb{P}\left(\frac{\log \max_{1 \leq k \leq n} X_{k}}{\log n} \leq \frac{1}{\rho} - y \right)~~\forall~ 0 < y < \frac{1}{\rho}.
\]

{\bf Step 4}.~~Under some reasonable conditions on $F(\cdot)$, in Theorem 2.5 we will obtain a result on convergence 
in distribution for $\log \max_{1 \leq k \leq n}X_{k}$, $n \geq 1$. 

{\bf Step 5}.~~Replacing $\max_{1 \leq k \leq n}X_{k}$ by $\max_{1 \leq k \leq n}|X_{k}|$ and following from Theorems 2.1-2.5, 
in Theorem 2.6 we will state a set of asymptotic properties of $\hat{\theta}_{n}$, $n \geq 1$. In particular, one of them asserts that 
\begin{equation}
\hat{\theta}_{n} \rightarrow_{\mathbb{P}} \theta
\end{equation}
if and only if 
\[
\lim_{x \rightarrow \infty} x^{r} \mathbb{P}(|X| > x) = \infty ~~\forall~r > \theta ~~\mbox{if}~~\theta < \infty,
\]
where ``$\rightarrow_{\mathbb{P}}$" stands for convergence in probability. If (1.2) holds for some $0 < \theta < \infty$, we
will see from Theorem 2.6 that
\[
\lim_{n \rightarrow \infty} \frac{\ln \mathbb{P}\left(\left| \hat{\theta}_{n} - \theta \right|  \geq \epsilon \right)}{\ln n} 
= - \frac{\epsilon}{\theta - \epsilon}~~\forall ~0 < \epsilon < \theta.
\]
This means that, under very reasonable conditions on $F(\cdot)$, $\hat{\theta}_{n}$ is not only a consistent estimator of $\theta$ 
but also possesses a very good convergence rate.  

The proofs of our main results will be provided in Section 3. As one can see from Section 3, the proofs of the main results are 
simple since only some basic results (such as the Borel-Cantelli lemma) in probability theory are used. We refer the reader to 
Chow and Teicher (1997) for any basic results in probability theory that are used in this paper.

In Section 4 hypothesis testing for the power of moments is conducted and, as an application of our main results, the formula 
for finding the p-value of the test is given. In addition, a theoretical application of our main results will be provided in 
Section 5 together with three illustrative examples.

\section{Statement of the main results}

Throughout, $X$ is a random variable with unknown distribution $F(x) = \mathbb{P}(X \leq x)$,
$-\infty < x < \infty$ and write
\[
\rho_{1} = \sup \left \{r \geq 0:~\lim_{x \rightarrow \infty} x^{r} \mathbb{P}(X > x) = 0 \right \}
~~\mbox{and}~~
\rho_{2} = \sup \left \{r \geq 0:~\liminf_{x \rightarrow \infty} x^{r} \mathbb{P}(X > x) = 0 \right \}.
\]
Clearly, just as $\theta$ as defined in Section 1 is a parameter of the distribution $F(\cdot)$ of the random variable $X$,
so are $\rho_{1}$ and $\rho_{2}$. These parameters satisfy
\[
0 \leq \rho_{1} \leq \rho_{2} \leq \infty.
\]
The main results of this paper are the following Theorems 2.1-2.6. 

\vskip 0.3cm

\begin{theorem}
Let $\{X_{n}; n \geq 1\}$ be a sequence of i.i.d. random variables drawn from the distribution function $F(\cdot)$ 
of the random variable $X$. Then 
\begin{equation}
\limsup_{n \rightarrow \infty} \frac{\log \max_{1 \leq k \leq n} X_{k}}{\log n} = \frac{1}{\rho_{1}}~~\mbox{a.s.}
\end{equation} 
and there exists an increasing positive integer sequence $\{l_{n}; n \geq 1 \}$ (which depends on the probability 
distribution of $X$ when $\rho_{1} < \infty$) such that
\begin{equation}
\lim_{n \rightarrow \infty} \frac{\log \max_{1 \leq k \leq l_{n}} X_{k}}{\log l_{n}} = \frac{1}{\rho_{1}}~~\mbox{a.s.}
\end{equation} 
\end{theorem}

\vskip 0.3cm

\begin{theorem}
Let $\{X_{n}; n \geq 1\}$ be a sequence of i.i.d. random variables drawn from the distribution function
$F(\cdot)$ of the random variable $X$. Then
\begin{equation}
\liminf_{n \rightarrow \infty} \frac{\log \max_{1 \leq k \leq n} X_{k}}{\log n} = \frac{1}{\rho_{2}}~~\mbox{a.s.}
\end{equation} 
and there exists an increasing positive integer sequence $\{m_{n}; n \geq 1 \}$ (which depends on the probability 
distribution of $X$ when $\rho_{2} > 0$) such that
\begin{equation}
\lim_{n \rightarrow \infty} \frac{\log \max_{1 \leq k \leq m_{n}} X_{k}}{\log m_{n}} = \frac{1}{\rho_{2}}~~\mbox{a.s.}
\end{equation} 
\end{theorem}

\vskip 0.3cm

\begin{remark}
We must point out that (2.2) and (2.4) are two interesting conclusions. To see this,
let $\{U_{n};~n \geq 1 \}$ be a sequence of independent random variables 
with
\[
\mathbb{P}\left(U_{n} = 1 \right) = \mathbb{P}\left(U_{n} = 3 \right) = \frac{1}{2n}
~~\mbox{and}~~\mathbb{P}\left(U_{n} = 2 \right) = 1 - \frac{1}{n}, ~~n \geq 1.
\]
Since 
\[
\sum_{n=1}^{\infty} \mathbb{P}\left(U_{n} = 3 \right) = 
\sum_{n=1}^{\infty} \mathbb{P}\left(U_{n} = 1 \right) = \sum_{n=1}^{\infty} \frac{1}{2n}
= \infty,
\]
it follows from the Borel-Cantelli lemma that
\[
\limsup_{n \rightarrow \infty} U_{n} = 3~~\mbox{a.s.}~~\mbox{and}~~
\liminf_{n \rightarrow \infty} U_{n} = 1~~\mbox{a.s.}
\]
However, for any sequences $\{l_{n};~n \geq 1\}$ and $\{m_{n};~n \geq 1 \}$
of increasing positive integers,
\[
\mbox{neither}~~\lim_{n \rightarrow \infty} U_{l_{n}} = 3~~\mbox{a.s.~~nor}~~
\lim_{n \rightarrow \infty} U_{m_{n}} = 1~~\mbox{a.s.~ holds.}
\]
\end{remark}

\vskip 0.3cm

\begin{remark}
For an observable random variable $X$, it is often the case that $\rho_{1} = \rho_{2}$. However, 
for any given constants $\rho_{1}$ and $\rho_{2}$ with $0 \leq \rho_{1} < \rho_{2} \leq \infty$,
one can construct a random variable $X$ such that
\[
\sup \left \{r \geq 0:~\limsup_{x \rightarrow \infty} x^{r} \mathbb{P}(X > x) = 0 \right \} = \rho_{1}
~~\mbox{and}~~\sup \left \{r \geq 0:~\liminf_{x \rightarrow \infty} x^{r} \mathbb{P}(X > x) = 0 \right \} = \rho_{2}.
\]
For example, if $0 < \rho_{1} < \rho_{2} < \infty$, a random variable $X$ can be constructed having probability
distribution given by
\[
\mathbb{P} \left(X = d_{n} \right) = \frac{c}{d_{n}^{\rho_{1}}}, ~~n \geq 1,
\]
where $d_{n} = 2^{\left(\rho_{2}/\rho_{1} \right)^{n}}$,  $n \geq 1$ and
\[
c = \left(\sum_{n = 1}^{\infty} \frac{1}{d_{n}^{\rho_{1}}} \right)^{-1} > 0.
\]
\end{remark}

\vskip 0.3cm

Combining Theorems 2.1 and 2.2, we establish a law of large numbers for $\log \max_{1 \leq k \leq n} X_{k}, n \geq 1$
as follows.

\vskip 0.3cm

\begin{theorem}
Let $\{X_{n}; n \geq 1\}$ be a sequence of i.i.d. random variables drawn from the distribution function $F(\cdot)$
of the random variable $X$ and let $\rho \in [0, \infty]$. Then the following four statements are equivalent:
\begin{equation}
\lim_{n \rightarrow \infty} \frac{\log \max_{1 \leq k \leq n} X_{k}}{\log n} = \frac{1}{\rho}~~\mbox{a.s.,}
\end{equation} 
\begin{equation}
\frac{\log \max_{1 \leq k \leq n} X_{k}}{\log n} \rightarrow_{\mathbb{P}} \frac{1}{\rho},
\end{equation}
\begin{equation}
\rho_{1} = \rho_{2} = \rho,
\end{equation}
\begin{equation}
\lim_{x \rightarrow \infty} x^{r} \mathbb{P}(X > x) =
\left \{
\begin{array}{ll}
0 & \forall~r < \rho ~\mbox{if}~ \rho > 0,\\
&\\
\infty & \forall~r > \rho ~\mbox{if}~ \rho < \infty.
\end{array}
\right.
\end{equation}
If $0 \leq \rho < \infty$, then anyone of (2.5)-(2.8) holds if and only if there exists a function 
$L(\cdot): (0, \infty) \rightarrow (0, \infty)$ such that
\begin{equation}
\mathbb{P}(X > x) \sim \frac{L(x)}{x^{\rho}}~~\mbox{as}~~ x \rightarrow \infty
~~\mbox{and}~~\lim_{x \rightarrow \infty} \frac{\ln L(x)}{\ln x} = 0.
\end{equation}
\end{theorem}

\vskip 0.3cm

The following result provides large deviation probabilities for $\log \max_{1 \leq k \leq n} X_{k}/\log n$, $n \geq 1$.

\begin{theorem}
Let $\{X_{n}; n \geq 1\}$ be a sequence of i.i.d. random variables drawn from the distribution function $F(\cdot)$
of the random variable $X$. If (2.6) holds for some $0 < \rho < \infty$, then
\begin{equation}
\lim_{n \rightarrow \infty} 
\frac{\ln \mathbb{P}\left(\frac{\log \max_{1 \leq k \leq n}X_{k}}{\log n} \geq \frac{1}{\rho} + x \right)}{\ln n}
= - \rho x ~~\forall~x > 0
\end{equation}
and
\begin{equation}
\lim_{n \rightarrow \infty} 
\frac{\log \left(-\ln \mathbb{P}
\left(\frac{\log \max_{1 \leq k \leq n}X_{k}}{\log n} \leq \frac{1}{\rho} - y \right) \right)}{\ln n}
= \rho y ~~\forall~0 < y < \frac{1}{\rho}.
\end{equation}
\end{theorem}

\vskip 0.3cm

\begin{remark}
If (2.6) holds for some $0 < \rho < \infty$, it then follows from (2.10) and (2.11) that
\[
\mathbb{P}\left(\frac{\log \max_{1 \leq k \leq n}X_{k}}{\log n} \geq \frac{1}{\rho} + x \right) 
= n^{-\rho x + o(1)} ~~\mbox{as}~ n \rightarrow \infty~~\forall~x > 0
\]	
and
\[
\mathbb{P}\left(\frac{\log \max_{1 \leq k \leq n}X_{k}}{\log n} \leq \frac{1}{\rho} - y \right) 
= \exp\left(-n^{\rho y + o(1)} \right) ~~\mbox{as}~~n \rightarrow \infty~~\forall~0 < y < \frac{1}{\rho}
\]
and hence 
\begin{equation}
\mathbb{P}\left(\left | \frac{\log \max_{1 \leq k \leq n}X_{k}}{\log n} - \frac{1}{\rho} \right| \geq \epsilon \right) 
= n^{-\rho \epsilon + o(1)} ~~\mbox{as}~~n \rightarrow \infty ~~\forall ~\epsilon > 0.
\end{equation}
\end{remark}

\vskip 0.3cm

The following result concerns convergence in distribution for $\log \max_{1 \leq k \leq n}X_{k}, ~n \geq 1$.

\vskip 0.3cm

\begin{theorem}
Let $\{X_{n}; n \geq 1\}$ be a sequence of i.i.d. random variables drawn from the distribution function $F(\cdot)$
of the random variable $X$. Suppose that there exist constants $0 < \rho < \infty$ and $-\infty < \tau < \infty$
and a monotone function $h(\cdot):~[0, \infty) \rightarrow (0, \infty)$ with $\lim_{x \rightarrow \infty}h(x^{2})/h(x) = 1$
such that
\begin{equation}
\mathbb{P}(X > x) \sim \frac{(\log x)^{\tau}h(x)}{x^{\rho}} ~~\mbox{as}~~x \rightarrow \infty.
\end{equation}
Then
\begin{equation}
\lim_{n \rightarrow \infty} \mathbb{P} 
\left(\log \max_{1 \leq k \leq n} X_{k} \leq \frac{\ln n + \tau \ln \ln n + \ln h(n) - \tau \ln \rho + x}{\rho} \right)
= \exp\left(-e^{-x} \right) ~~\forall~ - \infty < x < \infty.
\end{equation}
\end{theorem}

\vskip 0.3cm

We now return to the problem posed in Section 1. Note that, for $r > 0$
\[
\mbox{if}~~\lim_{x \rightarrow \infty} x^{r} \mathbb{P}(|X| > x) = 0
~~\mbox{then}~~\mathbb{E}|X|^{r_{1}} < \infty ~~\forall~0 \leq r_{1} < r
\]
and 
\[
\mbox{if}~~\mathbb{E}|X|^{r} < \infty
~~\mbox{then}~~\lim_{x \rightarrow \infty} x^{r_{1}} \mathbb{P}(|X| > x) = 0 ~~\forall~0 \leq r_{1} \leq r.
\]
We thus have that
\[
\sup\left \{r \geq 0:~\lim_{x \rightarrow \infty} x^{r} \mathbb{P}(|X| > x) = 0 \right \} 
= \sup\left \{r \geq 0:~\mathbb{E}|X|^{r} < \infty \right \} = \theta.
\]
Thus, by Theorems 2.1-2.5, some asymptotic properties of the point estimator $\hat{\theta}_{n}$ are provided 
in the following theorem.

\vskip 0.3cm

\begin{theorem}
Let $\{X_{n}; n \geq 1\}$ be a sequence of i.i.d. random variables drawn from the distribution function $F(\cdot)$
of the random variable $X$. Let
\[
\hat{\theta}_{n} = \frac{\log n}{\log \max_{1 \leq k \leq n} |X_{k}|}, ~n \geq 1.
\]
Then we have:

{\bf (i)}
\[
\liminf_{n \rightarrow \infty} \hat{\theta}_{n} = \theta = \sup \left \{r \geq 0:~ \mathbb{E}|X|^{r} < \infty \right \}~~\mbox{a.s.,}
\]
\[
\limsup_{n \rightarrow \infty}
\hat{\theta}_{n} = \sup \left\{r \geq 0:~ \liminf_{x \rightarrow \infty} x^{r} \mathbb{P}(|X| > x) = 0 \right \}~~\mbox{a.s.,}
\]
and the following three statements are equivalent:
\begin{equation}
\lim_{n \rightarrow \infty} \hat{\theta}_{n} = \theta~~\mbox{a.s.,}
\end{equation}
\begin{equation}
 \hat{\theta}_{n} \rightarrow_{\mathbb{P}} \theta,
\end{equation}
\begin{equation}
\lim_{x \rightarrow \infty} x^{r} \mathbb{P}(|X| > x) = \infty ~~\forall~r > \theta ~\mbox{if}~ \theta < \infty.
\end{equation}
If $0 \leq \theta < \infty$, then anyone of (2.15)-(2.17) holds if and only if there exists a function 
$L(\cdot): (0, \infty) \rightarrow (0, \infty)$ such that
\[
\mathbb{P}(|X| > x) \sim \frac{L(x)}{x^{\theta}}~~\mbox{as}~~ x \rightarrow \infty
~~\mbox{and}~~\lim_{x \rightarrow \infty} \frac{\ln L(x)}{\ln x} = 0.
\]

{\bf (ii)}~~If (2.16) holds for some $0 < \theta < \infty$, then
\[
\lim_{n \rightarrow \infty} 
\frac{\ln \mathbb{P}\left(\hat{\theta}_{n} \leq \theta - s \right)}{\ln n}
= - \frac{s}{\theta - s} ~~\forall~0 < s < \theta
\]
and
\[
\lim_{n \rightarrow \infty} 
\frac{\log \left(-\ln \mathbb{P}\left(\hat{\theta}_{n} \geq \theta + t \right) \right)}{\ln n}
= \frac{t}{\theta + t} ~~\forall~ t > 0
\]
and hence
\begin{equation}
\lim_{n \rightarrow \infty} 
\frac{\ln \mathbb{P} \left( \left |\hat{\theta}_{n} - \theta \right| \geq \epsilon \right)}{\ln n} 
= - \frac{\epsilon}{\theta - \epsilon}~~\forall~ 0 < \epsilon < \theta.
\end{equation}

{\bf (iii)}~~Suppose that there exist constants $0 < \theta < \infty$ and $-\infty < \tau < \infty$
and a monotone function $h(\cdot):~[0, \infty) \rightarrow (0, \infty)$ with $\lim_{x \rightarrow \infty}h(x^{2})/h(x) = 1$
such that
\[
\mathbb{P}(|X| > x) \sim \frac{(\log x)^{\tau}h(x)}{x^{\theta}} ~~\mbox{as}~~x \rightarrow \infty.
\]
Then
\[
\lim_{n \rightarrow \infty} \mathbb{P} 
\left(\hat{\theta}_{n} - \theta \leq \frac{-\theta \tau \ln \ln n - \theta \ln h(n) 
+ \theta \tau\ln \theta + \theta x}{\ln n} \right)
= 1 - \exp \left(-e^{x} \right)
~~\forall~ - \infty < x < \infty.
\]
\end{theorem}

\vskip 0.3cm

\begin{remark}
From Theorem 2.6, one can see that the point estimator $\hat{\theta}_{n}$ posseses some nice
asymptotic properties. In particular, it follows from (2.18) that
\[
\mathbb{P} \left(\left|\hat{\theta}_{n} - \theta \right | \geq \epsilon \right) 
= n^{-\epsilon/(\theta - \epsilon) + \mathfrak{o}(1)} 
~~\mbox{as}~~n \rightarrow \infty~~\forall~ 0 < \epsilon < \theta.
\]
Thus, under very reasonable conditions on $F(\cdot)$, $\hat{\theta}_{n}$ is a good candidate to be used for estimating $\theta$ 
since it is not only a consistent estimator of $\theta$ but also possesses a very good convergence rate.  
\end{remark}

\section{Proofs of the main results}

Let $\{A_{n};~n \geq 1 \}$ be a sequence of events. As usual the abbreviation $\left \{A_{n} ~\mbox{i.o.} \right \}$ 
stands for the event that the events $A_{n}$ occur infinitely often. That is,
\[
\left \{A_{n} ~\mbox{i.o.} \right \} = \left\{\mbox{events}~A_{n}~\mbox{occur infinitely often} \right \} 
= \bigcap_{n=1}^{\infty} \bigcup_{j=n}^{\infty} A_{j}.
\]
For events $A$ and $B$, we say $A = B$ a.s. if $\mathbb{P}(A \Delta B) = 0$ where 
$A \Delta B = (A \setminus B) \cup (B \setminus A)$. To prove Theorem 2.1, we use the following preliminary result
which can be found in Chandra (2012, Example 1.6.25 (a), p. 48).

\vskip 0.3cm

\begin{lemma}
Let $\left \{b_{n};~n \geq 1 \right \}$ be a  nondecreasing sequence of positive real numbers such that 
\[
\lim_{n \rightarrow \infty} b_{n} = \infty
\]
and let $\left \{V_{n}; ~n \geq 1 \right \}$ be a sequence of random variables. Then
\[
\left \{\max_{1 \leq k \leq n} V_{k} \geq b_{n} ~\mbox{i.o.} \right \} 
= \left \{ V_{n} \geq b_{n} ~\mbox{i.o.} \right \} ~~\mbox{a.s.}
\]
\end{lemma}

\vskip 0.3cm

{\it Proof of Theorem 2.1}~~{\bf Case I: $0 < \rho_{1} < \infty$.} For given $\epsilon > 0$, 
let $r(\epsilon) = \left(\frac{1}{\rho_{1}} + \epsilon \right)^{-1}$. Then 
\[
0 < r(\epsilon) < \rho_{1} = \sup \left \{r \geq 0:~\lim_{x \rightarrow \infty} x^{r} \mathbb{P}(X > x) = 0 \right \}
\]
and hence 
\begin{equation}
\sum_{n = 1}^{\infty} \mathbb{P} \left(X > n^{1/r(\epsilon)} \right) < \infty.
\end{equation}
By the Borel-Cantelli lemma, (3.1) implies that
\[
\mathbb{P} \left(X_{n} > n^{1/r(\epsilon)} ~\mbox{i.o.} \right) = 0.
\]
By Lemma 3.1, we have
\[
\left \{\frac{\log \max_{1 \leq k \leq n}X_{k}}{\log n} > \frac{1}{\rho_{1}} + \epsilon ~\mbox{i.o.} \right \}
= \left \{ \max_{1 \leq k \leq n} X_{k} > n^{1/r(\epsilon)} ~\mbox{i.o.} \right \} 
= \left \{X_{n} > n^{1/r(\epsilon)} ~\mbox{i.o.} \right \}~~\mbox{a.s.}
\]
and hence
\[
\mathbb{P} \left( \frac{\log \max_{1 \leq k \leq n}X_{k}}{\log n} > \frac{1}{\rho_{1}} + \epsilon ~\mbox{i.o.} \right) = 0.
\]
Thus 
\[
\limsup_{n \rightarrow \infty} \frac{\log \max_{1 \leq k \leq n}X_{k}}{\log n} \leq \frac{1}{\rho_{1}} + \epsilon ~~\mbox{a.s.}
\]
Letting $\epsilon \searrow 0$, we get
\begin{equation}
\limsup_{n \rightarrow \infty} \frac{\log \max_{1 \leq k \leq n}X_{k}}{\log n} \leq \frac{1}{\rho_{1}}~~\mbox{a.s.}
\end{equation}
By the definition of $\rho_{1}$, we have that
\[
\limsup_{x \rightarrow \infty} x^{r} \mathbb{P}(X > x) = \infty~~\forall~r > \rho_{1}
\]
which is equivalent to 
\[
\limsup_{x \rightarrow \infty} x \mathbb{P}\left(X > x^{(1/\rho_{1}) - \epsilon} \right)
= \infty~~\forall ~\epsilon > 0.
\]
Then, inductively, we can choose positive integers 
$l_{n}, n \geq 1$ such that
\[
1 = l_{1} < l_{2} < ... < l_{n} < ... ~~\mbox{and}~~l_{n}\mathbb{P}\left(X > l_{n}^{(1/\rho_{1}) - (1/n)} \right) \geq 2 \ln n, ~~n \geq 1.
\]
Note that, for any $0 \leq z \leq 1$, $1 - z \leq e^{-z}$. Thus, for all sufficiently large $n$, we have that
\[
\begin{array}{lll}
\mbox{$\displaystyle 
\mathbb{P} \left(\frac{\log \max_{1 \leq k \leq l_{n}} X_{k}}{\log l_{n}} \leq \frac{1}{\rho_{1}} - \frac{1}{n} \right)$}
&=& \mbox{$\displaystyle \mathbb{P} \left(\max_{1 \leq k \leq l_{n}}X_{k} \leq l_{n}^{(1/\rho_{1}) - (1/n)} \right)$}\\
&&\\
&=& 
\mbox{$\displaystyle \left(1 - \mathbb{P}\left(X > l_{n}^{(1/\rho_{1}) - (1/n)} \right) \right)^{l_{n}}$}\\
&&\\
&\leq &
\mbox{$\displaystyle \exp \left(- l_{n}\mathbb{P}\left(X > l_{n}^{(1/\rho_{1}) - (1/n)} \right) \right)$} \\
&&\\
&\leq&
\mbox{$\displaystyle \exp(-2 \ln n)$}\\
&&\\
&=&
\mbox{$\displaystyle n^{-2}$.}
\end{array}
\]
Since $\sum_{n=1}^{\infty} n^{-2} < \infty$, by the Borel-Cantelli lemma, we get that
\[
\mathbb{P} \left(\frac{\log \max_{1 \leq k \leq l_{n}} X_{k}}{\log l_{n}} \leq \frac{1}{\rho_{1}} - \frac{1}{n}~\mbox{i.o.} \right) = 0
\]
which ensures that
\begin{equation}
\liminf_{n \rightarrow \infty} \frac{\log \max_{1 \leq k \leq l_{n}} X_{k}}{\log l_{n}} \geq \frac{1}{\rho_{1}}~~\mbox{a.s.}
\end{equation}
Clearly, (2.1) and (2.2) follow from (3.2) and (3.3).

{\bf Case II: ~$\rho_{1} = \infty$.}~~Using the same argument used in the first half of the proof for Case I, we get that
\[
\limsup_{n \rightarrow \infty} \frac{\log \max_{1 \leq k \leq n} X_{k}}{\log n} \leq \epsilon~~\mbox{a.s.} ~~\forall~\epsilon > 0
\]
and hence
\begin{equation}
\limsup_{n \rightarrow \infty} \frac{\log \max_{1 \leq k \leq n} X_{k}}{\log n} \leq 0~~\mbox{a.s.} 
\end{equation}
Note that
\[
0 \leq \frac{\log \max_{1 \leq k \leq n} X_{k}}{\log n}~~\forall~n \geq 1.
\]
We thus have that
\begin{equation}
\liminf_{n \rightarrow \infty} \frac{\log \max_{1 \leq k \leq n} X_{k}}{\log n} \geq 0 ~~\mbox{a.s.}
\end{equation}
It thus follows from (3.4) and (3.5) that
\[
\lim_{n \rightarrow \infty} \frac{\log \max_{1 \leq k \leq n} X_{k}}{\log n} = 0 ~~\mbox{a.s.}
\]
proving (2.1) and (2.2) (with $l_{n} = n$, $n \geq 1$).

{\bf Case III: ~$\rho_{1} = 0$.}~~By the definition of $\rho_{1}$, we have that
\[
\limsup_{x \rightarrow \infty} x^{r} \mathbb{P}(X > x) = \infty~~\forall~r > 0
\]
which is equivalent to 
\[
\limsup_{x \rightarrow \infty} x \mathbb{P}\left(X > x^{r} \right)
= \infty~~\forall ~r > 0.
\]
Then, inductively, we can choose positive integers 
$l_{n}, n \geq 1$ such that
\[
1 = l_{1} < l_{2} < ... < l_{n} < ... ~~\mbox{and}~~l_{n}\mathbb{P}\left(X > l_{n}^{n} \right) \geq 2 \ln n, ~~n \geq 1.
\]
Thus, for all sufficiently large $n$, we have by the same argument as in Case I that
\[
\mathbb{P} \left(\frac{\log \max_{1 \leq k \leq l_{n}}X_{k}}{\log l_{n}} \leq n \right )
\leq n^{-2}
\]
and hence by the Borel-Cantelli lemma
\[
\mathbb{P} \left(\frac{\log \max_{1 \leq k \leq l_{n}} X_{k}}{\log l_{n}} \leq n ~\mbox{i.o.} \right) = 0
\]
which ensures that
\[
\lim_{n \rightarrow \infty} \frac{\log \max_{1 \leq k \leq l_{n}} X_{k}}{\log l_{n}} = \infty~~\mbox{a.s.}
\]
Thus (2.1) and (2.2) hold. This completes the proof of Theorem 2.1. ~$\Box$

\vskip 0.3cm

{\it Proof of Theorem 2.2}~~~~{\bf Case I: $0 < \rho_{2} < \infty$.} For given $ \rho_{2} < r < \infty$, 
let $r_{1} = \left(r + \rho_{2} \right)/2$ and $\tau = 1 - (r_{1}/r)$. Then $\rho_{2} < r_{1} < r < \infty$
and $\tau > 0$. By the definition of $\rho_{2}$, we have that
\[
\lim_{x \rightarrow \infty} x^{r_{1}} \mathbb{P}(X > x) = \infty
\]
and hence for all sufficiently large $x$,
\[
\mathbb{P}(X > x) \geq x^{-r_{1}}.
\]
Thus, for all sufficiently large $n$,
\[
n \mathbb{P}\left(X > n^{1/r} \right) \geq n \left(n^{1/r} \right)^{-r_{1}} = n^{1 - (r_{1}/r)} = n^{\tau}
\]
and hence
\[
\mathbb{P} \left(\max_{1 \leq k \leq n}X_{k} \leq n^{1/r} \right) = \left(1 - \mathbb{P}\left(X > n^{1/r} \right) \right)^{n}
\leq e^{-n \mathbb{P}\left(X > n^{1/r} \right)} \leq e^{-n^{\tau}}.
\]
Since
\[
\sum_{n=1}^{\infty} e^{-n^{\tau}} < \infty,
\]
by the Borel-Cantelli lemma, we have that
\[
\mathbb{P} \left(\max_{1 \leq k \leq n} X_{k} \leq n^{1/r} ~\mbox{i.o.} \right) = 0
\]
which implies that
\[
\liminf_{n \rightarrow \infty} \frac{\log \max_{1 \leq k \leq n}X_{k}}{\log n} \geq 1/r ~~\mbox{a.s.}
\]
Letting $r \searrow \rho_{2}$, we get
\begin{equation}
\liminf_{n \rightarrow \infty} \frac{\log \max_{1 \leq k \leq n}X_{k}}{\log n} \geq \frac{1}{\rho_{2}}~~\mbox{a.s.}
\end{equation}
Again, by the definition of $\rho_{2}$, we have that
\[
\liminf_{x \rightarrow \infty} x^{r} \mathbb{P}(X > x) = 0~~\forall~r < \rho_{2}
\]
which is equivalent to 
\[
\liminf_{x \rightarrow \infty} x \mathbb{P}\left(X > x^{(1/\rho_{2}) + \epsilon} \right)
= 0~~\forall ~\epsilon > 0.
\]
Then, inductively, we can choose positive integers 
$m_{n}, n \geq 1$ such that
\[
1 = m_{1} < m_{2} < ... < m_{n} < ... ~~\mbox{and}~~m_{n}\mathbb{P}\left(X > m_{n}^{(1/\rho_{2}) + (1/n)} \right) \leq n^{-2}, ~~n \geq 1.
\]
Then we have that
\[
\sum_{n=1}^{\infty} \mathbb{P}\left(\max_{1 \leq k \leq m_{n}}X_{k} >  m_{n}^{(1/\rho_{2}) + (1/n)} \right)
\leq \sum_{n=1}^{\infty} m_{n}\mathbb{P}\left(X > m_{n}^{(1/\rho_{2}) + (1/n)} \right) \leq \sum_{n=1}^{\infty} n^{-2} < \infty.
\]
Thus, by the Borel-Cantelli lemma, we get that
\[
\mathbb{P} \left(\frac{\log \max_{1 \leq k \leq m_{n}} X_{k}}{\log m_{n}} > \frac{1}{\rho_{2}} + \frac{1}{n}~\mbox{i.o.} \right) = 0
\]
which ensures that
\begin{equation}
\limsup_{n \rightarrow \infty} \frac{\log \max_{1 \leq k \leq m_{n}} X_{k}}{\log m_{n}} \leq \frac{1}{\rho_{2}}~~\mbox{a.s.}
\end{equation}
Clearly, (2.3) and (2.4) follow from (3.6) and (3.7).

{\bf Case II: ~$\rho_{2} = \infty$.}~~By the definition of $\rho_{2}$, we have that
\[
\liminf_{x \rightarrow \infty} x^{r} \mathbb{P}(X > x) = 0~~\forall~r > 0
\]
which is equivalent to 
\[
\liminf_{x \rightarrow \infty} x \mathbb{P}\left(X > x^{r} \right)
= 0~~\forall ~r > 0.
\]
Then, inductively, we can choose positive integers 
$m_{n}, n \geq 1$ such that
\[
1 = m_{1} < m_{2} < ... < m_{n} < ... ~~\mbox{and}~~m_{n}\mathbb{P}\left(X > m_{n}^{1/n} \right) \leq n^{-2}, ~~n \geq 1.
\]
Thus
\[
\sum_{n=1}^{\infty} \mathbb{P} \left(\max_{1 \leq k \leq m_{n}}X_{k} > m_{n}^{1/n} \right) 
\leq \sum_{n=1}^{\infty} m_{n}\mathbb{P}\left(X > m_{n}^{1/n} \right) \leq \sum_{n=1}^{\infty} n^{-2} < \infty
\]
and hence the Borel-Cantelli lemma
\[
\mathbb{P} \left( \max_{1 \leq k \leq m_{n}} X_{k} > m_{n}^{1/n} ~\mbox{i.o.} \right) = 0
\]
which ensures that
\begin{equation}
\limsup_{n \rightarrow \infty} \frac{\log \max_{1 \leq k \leq m_{n}} X_{k}}{\log m_{n}} \leq 0~~\mbox{a.s.}
\end{equation}
It is clear that
\begin{equation}
\liminf_{n \rightarrow \infty} \frac{\log \max_{1 \leq k \leq n} X_{k}}{\log n} \geq 0~~\mbox{a.s.}
\end{equation}
It thus follows from (3.8) and (3.9) that
\[
\liminf_{n \rightarrow \infty} \frac{\log \max_{1 \leq k \leq n} X_{k}}{\log n} = 0~~\mbox{a.s.}
~~\mbox{and}~~\lim_{n \rightarrow \infty} \frac{\log \max_{1 \leq k \leq m_{n}} X_{k}}{\log m_{n}} = 0~~\mbox{a.s.;}
\]
i.e., (2.3) and (2.4) hold. 

{\bf Case III: ~$\rho_{2} = 0$.}~~Using the same argument used in the first half of the proof for Case I, we get that
\[
\liminf_{n \rightarrow \infty} \frac{\log \max_{1 \leq k \leq n} X_{k}}{\log n} \geq \frac{1}{r}~~\mbox{a.s.} ~~\forall~r > 0.
\]
Letting $r \searrow 0$, we get that
\[
\liminf_{n \rightarrow \infty} \frac{\log \max_{1 \leq k \leq n} X_{k}}{\log n} = \infty~~\mbox{a.s.} 
\]
Thus
\[
\lim_{n \rightarrow \infty} \frac{\log \max_{1 \leq k \leq n} X_{k}}{\log n} = \infty ~~\mbox{a.s.}
\]
and hence (2.3) and (2.4) hold with $m_{n} = n$, $n \geq 1$. ~$\Box$

\vskip 0.3cm

{\it Proof of Theorem 2.3}~~It follows from Theorems 2.1 and 2.2 that
\[
(2.5) \Longleftrightarrow (2.7) \Longleftrightarrow (2.8).
\]
Since (2.6) follows from (2.5), we only need to show that (2.6) implies (2.8). It follows from (2.6) that 
\begin{equation}
\lim_{n \rightarrow \infty} \mathbb{P} \left( \frac{\log \max_{1 \leq k \leq n} X_{k}}{\log n} \leq \frac{1}{r} \right)
= \left \{
\begin{array}{ll}
1 & \forall~r < \rho~ \mbox{if}~ \rho > 0,\\
&\\
0 & \forall~r > \rho ~\mbox{if}~ \rho < \infty.
\end{array}
\right.
\end{equation}
Since, for $n \geq 3$
\[
\mathbb{P} \left( \frac{\log \max_{1 \leq k \leq n} X_{k}}{\log n} \leq \frac{1}{r} \right)
= \mathbb{P} \left( \max_{1 \leq k \leq n} X_{k} \leq n^{1/r} \right)
= \left(1 - \mathbb{P}\left(X > n^{1/r} \right) \right)^{n} = e^{n \ln \left(1 - \mathbb{P}\left(X > n^{1/r} \right) \right)}
\]
and
\[
n \ln \left(1 - \mathbb{P}\left(X > n^{1/r} \right)  \right) \sim - n \mathbb{P}\left(X > n^{1/r} \right)~~\mbox{as}~~n \rightarrow \infty,
\]
it follows from (3.10) that
\[
\lim_{n \rightarrow \infty} n \mathbb{P}\left(X > n^{1/r} \right)
= \left \{
\begin{array}{ll}
0 & \forall~r < \rho~ \mbox{if}~ \rho > 0,\\
&\\
\infty & \forall~r > \rho~ \mbox{if}~ \rho < \infty
\end{array}
\right.
\]
which is equivalent to (2.8). 

For $0 \leq \rho < \infty$, note that
\[
\mathbb{P} (X > x) = x^{-\rho} \left(x^{\rho}\mathbb{P}(X > x) \right) 
= e^{-\rho \ln x + \ln \left(x^{\rho} \mathbb{P}(X > x) \right)}~~\forall ~x > 0.
\]
We thus see that, if $0 \leq \rho < \infty$, then (2.8) is equivalent to 
\[
\lim_{x \rightarrow \infty} \frac{\ln \left(x^{\rho} \mathbb{P}(X > x) \right)}{\log x}
= 0.
\]
(We leave it to the reader to work out the details of the proof.) We thus see that (2.8) implies
(2.9) with $L(x) = \ln \left(x^{\rho} \mathbb{P}(X > x) \right)$, $x > 0$. It is easy to verify 
that (2.8) follows from (2.9). This completes the proof of Theorem 2.3. ~$\Box$

\vskip 0.3cm

{\it Proof of Theorem 2.4}~~Since (2.6) holds for some $0 < \rho < \infty$, it follows from the proof 
of Theorem 2.3 that the function $L(x) = x^{\rho}\mathbb{P}(X > x)$, $x > 0$ satisfies 
\[
\lim_{x \rightarrow \infty} \frac{\ln L(x)}{\ln x} = 0.
\]
Thus, for fixed $x > 0$ and $0 < y < 1/\rho$, we have that, as $n \rightarrow \infty$,
\[
\begin{array}{lll}
\mbox{$\displaystyle
\mathbb{P} \left(\frac{\log \max_{1 \leq k \leq n} X_{k}}{\log n} \geq \frac{1}{\rho} + x \right) $}
&=& \mbox{$\displaystyle 
1 - \left(1 - \mathbb{P}\left(X \geq n^{(1/\rho) + x} \right) \right)^{n} $}\\
&&\\
&=& \mbox{$\displaystyle 
1 - e^{n \ln \left(1 - \mathbb{P}\left(X \geq n^{(1/\rho) + x} \right) \right)}$}\\
&&\\
&\sim &  \mbox{$\displaystyle -n \ln \left(1 - \mathbb{P}\left(X \geq n^{(1/\rho) + x} \right) \right)$}\\
&&\\
&\sim & \mbox{$\displaystyle n \mathbb{P}\left(X \geq n^{(1/\rho) + x} \right)$}\\
&&\\
&\sim& \mbox{$\displaystyle 
n \frac{L\left(n^{(1/\rho) + x}\right)}{n^{\rho\left(\left(\frac{1}{\rho}\right) + x \right)}}$}\\
&&\\
&=& \mbox{$\displaystyle
n^{-\rho x} L\left(n^{(1/\rho) + x}\right)$}
\end{array}
\]
and 
\[
\begin{array}{lll}
\mbox{$\displaystyle
\mathbb{P} \left(\frac{\log \max_{1 \leq k \leq n} X_{k}}{\log n} \leq \frac{1}{\rho} - y \right) $}
&=& \mbox{$\displaystyle 
\left(1 - \mathbb{P}\left(X \geq n^{(1/\rho) - y} \right) \right)^{n}$}\\
&&\\
&=& \mbox{$\displaystyle 
e^{n \ln \left(1 -\mathbb{P}\left(X \geq n^{(1/\rho) - y} \right) \right)}$}\\
&&\\
& \sim & \mbox{$\displaystyle
e^{-n \mathbb{P}\left(X \geq n^{(1/\rho) - y} \right)}$} \\
&&\\
&\sim& \mbox{$\displaystyle 
e^{-n^{\rho y} L\left(X > n^{(1/\rho) - y} \right)}$.}
\end{array}
\]
We thus have that
\[
\lim_{n \rightarrow \infty} 
\frac{\ln \mathbb{P} \left(\frac{\log \max_{1 \leq k \leq n} X_{k}}{\log n} \geq \frac{1}{\rho} + x \right)}{\ln n } 
= -\rho x + \lim_{n \rightarrow \infty} \frac{\ln L\left(n^{(1/\rho) + x} \right)}{\ln n} = - \rho x
\]
and
\[
\lim_{n \rightarrow \infty} 
\frac{\log \left(-\ln \mathbb{P}\left(\frac{\log \max_{1 \leq k \leq n}X_{k}}{\log n} \leq \frac{1}{\rho} - y \right) \right)}{\ln n}
= \rho y + \lim_{n \rightarrow \infty} \frac{L\left(X > n^{(1/\rho) - y} \right)}{\ln n}
= \rho y;
\]
i.e., (2.10) and (2.11) hold. ~$\Box$

\vskip 0.3cm

{\it Proof of Theorem 2.5}~~For fixed $x \in (-\infty, \infty)$, write
\[
a_{n}(x) = \frac{\ln n + \tau \ln \ln n + \ln h(n) - \tau \ln \rho + x}{\rho}~~\mbox{and}~~b_{n}(x) = e^{a_{n}(x)}, ~n \geq 2.
\]
Then
\[
b_{n}(x) = n^{1/\rho} (\ln n)^{\tau/\rho}(h(n))^{1/\rho} \rho^{-\tau/\rho} e^{x/\rho}, ~n \geq 2.
\]
Since $h(\cdot):~[0, \infty) \rightarrow (0, \infty)$ is a monotone function with $\lim_{x \rightarrow \infty}h(x^{2})/h(x) = 1$,
$h(\cdot)$ is a slowly varying function such that $\lim_{x \rightarrow \infty} h(x^{r})/h(x) = 1$ $\forall~ r > 0$
and hence
\[
h\left(b_{n}(x) \right) \sim h(n) ~~\mbox{as}~~n \rightarrow \infty.
\]
Clearly,
\[
\left(\ln b_{n}(x) \right)^{\tau} \sim \rho^{-\tau}(\ln n)^{\tau}~~\mbox{as}~~n \rightarrow \infty.
\]
It thus follows from (2.13) that, as $n \rightarrow \infty$, 
\[
\begin{array}{lll}
\mbox{$\displaystyle
n \ln \left(1 - \mathbb{P}\left(X > b_{n}(x) \right) \right)$}
& \sim & 
\mbox{$\displaystyle 
-n \mathbb{P}\left(X > b_{n}(x) \right)$} \\
&&\\
& \sim & 
\mbox{$\displaystyle 
-n \times \frac{\left(\ln\left(b_{n}(x) \right) \right)^{\tau} h\left(b_{n}(x) \right)}{\left(b_{n}(x) \right)^{\rho}}$}\\
&&\\
& \sim &
\mbox{$\displaystyle 
-n \times \frac{\rho^{-\tau}(\ln n)^{\tau} h(n)}{n (\ln n)^{\tau} h(n) \rho^{-\tau} e^{x}}$}\\
&&\\
&=& 
\mbox{$\displaystyle -e^{-x}$}
\end{array} 
\]
so that
\[
\begin{array}{lll}
\mbox{$\displaystyle
\lim_{n \rightarrow \infty} \mathbb{P} \left(\log \max_{1 \leq k \leq n} X_{k} \leq a_{n}(x) \right)$}
&=& 
\mbox{$\displaystyle
\lim_{n \rightarrow \infty} \mathbb{P} \left(\max_{1 \leq k \leq n} X_{k} \leq b_{n}(x) \right)$}\\
&&\\
&=& \mbox{$\displaystyle
\lim_{n \rightarrow \infty} \left(1 - \mathbb{P} \left(X > b_{n}(x) \right) \right)^{n}$}\\
&&\\
&=&
\mbox{$\displaystyle
\lim_{n \rightarrow \infty} e^{n \ln \left(1 - \mathbb{P} \left(X > b_{n}(x) \right) \right)}$}\\
&&\\
&=& 
\mbox{$\displaystyle 
\exp \left(-e^{-x} \right)$;}
\end{array}
\]
i.e., (2.14) holds. ~$\Box$

\vskip 0.3cm

{\it Proof of Theorem 2.6}~~Since $\displaystyle \hat{\theta}_{n} 
= \frac{\log n}{\log \max_{1 \leq k \leq n} \left|X_{k} \right|}$, $n \geq 1$, 
Theorem 2.6 (i) follows immediately from Theorems 2.1-2.3. 

Since 
\[
\mathbb{P}\left(\hat{\theta}_{n} \leq \theta - s \right) 
= \mathbb{P}\left(\frac{\log \max_{1 \leq k \leq n}\left|X_{k} \right|}{\log n} \geq \frac{1}{\theta} + \frac{s}{\theta(\theta - s)} \right) 
~~\forall~ 0 < s < \theta
and 
\]
\[
\mathbb{P}\left(\hat{\theta}_{n} \geq \theta +  t \right) 
= \mathbb{P}\left(\frac{\log \max_{1 \leq k \leq n}\left|X_{k} \right|}{\log n} \leq \frac{1}{\theta} - \frac{t}{\theta(\theta + t)} \right) 
~~\forall~ t > 0,
\]
Theorem 2.6 (ii) follows from Theorem 2.4.

Under the conditions of Theorem 2.6 (iii), by Theorem 2.5 we have that
\[
\lim_{n \rightarrow \infty} \mathbb{P} 
\left(\log \max_{1 \leq k \leq n} \left|X_{k} \right| \leq \frac{\ln n + \tau \ln \ln n + \ln h(n) - \tau \ln \theta + x}{\theta} \right)
= \exp\left(-e^{-x} \right) ~~\forall~ - \infty < x < \infty
\]
and hence 
\begin{equation}
\lim_{n \rightarrow \infty} \mathbb{P} 
\left(\hat{\theta}_{n} \geq \frac{\theta \ln n}{\ln n + \tau \ln \ln n + \ln h(n) - \tau \ln \theta + x} \right)
= \exp\left(-e^{-x} \right) ~~\forall~ - \infty < x < \infty.
\end{equation}
Since $h(z):~[0, \infty) \rightarrow (0, \infty)$ is a monotone function with $\lim_{z \rightarrow \infty}h(z^{2})/h(z) = 1$,
$h\left(e^{z} \right)$ is a slowly varying function and hence
\[
z^{-1} \leq h\left(e^{z} \right) \leq z ~~\mbox{for all sufficiently large}~z;
\]
i.e., 
\[
(\ln z)^{-1} \leq h(z) \leq \ln z ~~\mbox{for all sufficiently large}~z.
\]
Thus
\[
- \ln \ln z \leq \ln h(z) \leq \ln \ln z ~~~~\mbox{for all sufficiently large}~z.
\]
Thus, for fixed $x$, we have that
\begin{equation}
\begin{array}{ll}
& \mbox{$\displaystyle 
\frac{\theta \ln n}{\ln n + \tau \ln \ln n + \ln h(n) - \tau \ln \theta + x}$}\\
&\\
& \mbox{$\displaystyle 
= \frac{\theta}{1 + \frac{\tau \ln \ln n + \ln h(n) - \tau \ln \theta + x}{\ln n}}$}\\
&\\
&\mbox{$\displaystyle 
=\theta \left( 1 - \frac{\tau \ln \ln n + \ln h(n) - \tau \ln \theta + x}{\ln n}
+ \mathcal{O} \left( \left(\frac{\ln \ln n}{\ln n} \right)^{2} \right)\right)$}\\
&\\
& \mbox{$\displaystyle 
= \theta + \frac{-\theta \tau \ln \ln n - \theta \ln h(n) + \theta \tau \ln \theta + \theta (-x)}{\ln n} 
+ o\left(\frac{1}{(\ln n)^{3/2}} \right)$.}
\end{array}
\end{equation}
It now follows from (3.11) and (3.12) that
\[
\lim_{n \rightarrow \infty} \mathbb{P} 
\left(\hat{\theta}_{n} - \theta \geq \frac{-\theta \tau \ln \ln n - \theta \ln h(n) + \theta \tau\ln \theta + \theta (-x)}{\ln n} \right)
= \exp \left(-e^{-x} \right)
~~\forall~ - \infty < x < \infty
\]
and hence
\[
\begin{array}{ll}
& \mbox{$\displaystyle
\lim_{n \rightarrow \infty} \mathbb{P} 
\left(\hat{\theta}_{n} - \theta \leq \frac{-\theta \tau \ln \ln n - \theta \ln h(n) + \theta \tau\ln \theta + \theta x}{\ln n} \right)$}\\
&\\
& \mbox{$\displaystyle
= 1 - \lim_{n \rightarrow \infty} \mathbb{P} 
\left(\hat{\theta}_{n} - \theta > \frac{-\theta \tau \ln \ln n - \theta \ln h(n) + \theta \tau\ln \theta + \theta (-(-x))}{\ln n} \right)$}\\
&\\
& \mbox{$\displaystyle 
= 1 - \exp \left(-e^{x} \right)
~~\forall~ - \infty < x < \infty$.}
\end{array}
\]
This proves Theorem 2.6 (iii). ~$\Box$

\section{Hypothesis testing for the power of moments}

We now return to the statistical problem addressed in Section 1. Let $X_{1}, X_{2}, ..., X_{n}$ be a random sample of size $n$ 
drawn from an observable random variable $X$ with unknown distribution function $F(x) = \mathbb{P}(X \leq x),
- \infty < x < \infty$. Let $\theta$ be the power of moments of the random variable $X$. Since, under very reasonable conditions 
on $F(\cdot)$, $\hat{\theta}_{n}$ is not only a consistent estimator of $\theta$ but also possesses a very good convergence rate,
we use $\hat{\theta}_{n}$ to estimate $\theta$. Let $\theta_{0} \in (0, \infty)$ be a specific value. In order to determine that
$\theta$ is greater than $\theta_{0}$, we conduct the following test of hypothesis for $\theta$:
\begin{equation}
H_{0}:~~\theta \leq \theta_{0} ~~\mbox{versus}~~H_{1}:~~\theta > \theta_{0}
\end{equation}
and use $\hat{\theta}_{n}$ to test (4.1). 

Let $\theta_{1}$ be the observed value of $\hat{\theta}_{n}$ based on an obtained data set. Then, for testing (4.1), under very reasonable 
conditions on $F(\cdot)$, it follows from Theorem 2.6 (ii) that
\begin{equation}
\begin{array}{lll}
\mbox{p-value} &=& 
\mbox{$\displaystyle \mathbb{P}\left(\hat{\theta}_{n} > \theta_{1} \left | \theta = \theta_{0} \right. \right)$}\\
&&\\
&=& \mbox{$\displaystyle 
\left \{
\begin{array}{ll}
\mathbb{P} \left(\hat{\theta}_{n} > \theta_{0} + \left(\theta_{1} - \theta_{0} \right) \left | \theta = \theta_{0} \right. \right) 
& \mbox{if}~\theta_{1} > \theta_{0},\\
&\\
1 - \mathbb{P} \left(\hat{\theta}_{n} \leq \theta_{0} - \left(\theta_{0} - \theta_{1} \right) \left | \theta = \theta_{0} \right. \right) 
& \mbox{if}~0 < \theta_{1} < \theta_{0}
\end{array}
\right. $}\\
&&\\
&=& \mbox{$\displaystyle 
\left \{
\begin{array}{ll}
\exp \left(- n^{\left(1 + o(1) \right) \frac{\theta_{1} - \theta_{0}}{\theta_{1}}} \right)
& \mbox{if}~\theta_{1} > \theta_{0},\\
&\\
1 - n^{-\left(1 + o(1) \right) \frac{\theta_{0} - \theta_{1}}{\theta_{1}}}
& \mbox{if}~0 < \theta_{1} < \theta_{0}.
\end{array}
\right. $}\\
\end{array}
\end{equation}
Let $\alpha$ be a given level of significance. If the calculated p-value is greater than $\alpha$, we then fail to reject the
null hypothesis $H_{0}:~ \theta \leq \theta_{1}$ at the $\alpha$ level of significance. Otherwise, there is sufficient evidence 
(at the $\alpha$ level of significance) to conclude that the alternative hypothesis $H_{1}:~\theta > \theta_{0}$ is true. 

Although the formula (4.2) can be used to calculate the p-value approximately for testing (4.1), it does not provide for us such a
formula for the case $\theta_{1} = \theta_{0}$. The following example shows us how the p-value can be found for the case
$\theta_{1} = \theta_{0}$.   

\vskip 0.3cm

\begin{example}
Let $X_{1}, X_{2}, ..., X_{n}$ be a random sample of size $n$ drawn from a population random variable $X$ such that
\[
\mathbb{P}(|X| > x) \sim \frac{c(\log x)^{\tau}}{x^{\theta}} ~~\mbox{as}~~x \rightarrow \infty,
\]
where $0 < \theta < \infty$, $0 < c < \infty$, and $-\infty < \tau < \infty$ are constants. 
For the case $\theta_{1} = \theta_{0}$, we have that, for all sufficiently large $n$
\[
\begin{array}{lll}
\mbox{\rm p-value} &=& 
\mbox{$\displaystyle \mathbb{P}\left(\hat{\theta}_{n} > \theta_{0} \left | \theta = \theta_{0} \right. \right)$}\\
&&\\
&=& \mbox{$\displaystyle 
\mathbb{P} \left(\max_{1 \leq k \leq n}|X_{k}| < \left. n^{1/\theta_{0}} \right | \theta = \theta_{0} \right)$}\\
&&\\
&=&  \mbox{$\displaystyle 
\left(1 - \left(1 + o(1) \right) \frac{\left(c/\theta_{0}^{\tau}\right)(\log n)^{\tau}}{n} \right)^{n}$}\\
&&\\
&\rightarrow & \mbox{$\displaystyle
\left \{
\begin{array}{ll}
0 & \mbox{if}~\tau > 0,\\
&\\
e^{-c} & \mbox{if}~\tau = 0,\\
&\\
1 & \mbox{if}~ \tau < 0.
\end{array}
\right. 
$}
\end{array}
\]

\end{example}

\section{A theoretical application of the main results}

Let $\{X_{n}; n \geq 1\}$ be a sequence of i.i.d. random variables drawn from the distribution function
$F(\cdot)$ of the random variable $X$. Then $\left\{\max_{1 \leq k \leq n} X_{k};~n \geq 1 \right \}$ is
called the {\it maxima sequence} associated with $\{X_{n}; n \geq 1\}$. Thus we see that our main results  
are actually stability theorems for the maxima sequence. The stability properties for the maxima sequence,
which are useful in many practical situations where we are interested in extreme behaviour rather than average 
behaviour, have been studied by Gnedenko (1943), Barndorff-Nielsen (1963), Tomkins (1986), and many other authors.

The following classical and well-known stability theorem is due to Barndorff-Nielsen (1963).

\vskip 0.3cm

\noindent {\bf Barndorff-Nielsen Stability Theorem}~~{\it Let $\{X_{n}; n \geq 1\}$ be a sequence of i.i.d. 
random variables drawn from the distribution function $F(\cdot)$ of the random variable $X$ with
\[
\sup\left\{x:~F(x) < 1 \right \} = \infty.
\]
Then there exists a sequence $\{a_{n};~n \geq 1 \}$ of real numbers such that
\begin{equation}
\lim_{n \rightarrow \infty} \frac{\max_{1 \leq k \leq n} X_{k}}{a_{n}} = 1~~\mbox{a.s.}
\end{equation}
if and only if
\begin{equation}
\int_{0}^{\infty} \frac{dF(x)}{1 - F(\delta x)} < \infty~~\forall~0 < \delta < 1.
\end{equation}
In either case, the sequence $\{a_{n};~n \geq 1\}$ may be assumed to be
\[
\mu_{n} = \inf\left \{x > 0:~F(x) \geq 1 - \frac{1}{n} \right \}, ~~n \geq 1.
\]
}

\vskip 0.2cm

Since it usually can be very complicated to check whether the integral in (5.2)  
is convergent or divergent and to find $\left\{\mu_{n}; n \geq 1 \right \}$, it is natural for us 
to seek an easy approach to see whether (5.1) holds and if so, to find $\{a_{n};~n \geq 1\}$ 
easily and quickly. As an application of our Theorem 2.3, in this section 
we will provide such a powerful method; see Theorem 5.2 below.

First, our main results will be used to establish the following stability theorem for the maxima sequence. 

\vskip 0.2cm

\begin{theorem}
Let $0 < \lambda < \infty$ and let $\beta \in [0, \infty]$. Let $\{X_{n}; n \geq 1\}$ be a sequence of i.i.d. 
random variables drawn from the distribution function $F(\cdot)$ of the random variable $X$. Write
\[
\beta_{1} = \sup \left \{r \geq 0:~\lim_{x \rightarrow \infty} e^{r x^{1/\lambda}} \mathbb{P}(X > x) = 0 \right \}~~\mbox{and}~~
\beta_{2} = \sup \left \{r \geq 0:~\liminf_{x \rightarrow \infty} e^{r x^{1/\lambda}} \mathbb{P}(X > x) = 0 \right \}.
\]
Then we have
\begin{equation}
\limsup_{n \rightarrow \infty} \frac{\max_{1 \leq k \leq n}X_{k}}{(\log n)^{\lambda}} = \beta_{1}^{-\lambda}~~\mbox{a.s. ~and}~~
\liminf_{n \rightarrow \infty} \frac{\max_{1 \leq k \leq n}X_{k}}{(\log n)^{\lambda}} = \beta_{2}^{-\lambda}~~\mbox{a.s.}
\end{equation}
and the following four statements are equivalent:
\begin{equation}
\lim_{n \rightarrow \infty} \frac{\max_{1 \leq k \leq n}X_{k}}{(\log n)^{\lambda}} = \beta^{-\lambda}~~\mbox{a.s.,}
\end{equation}
\begin{equation}
\frac{\max_{1 \leq k \leq n}X_{k}}{(\log n)^{\lambda}} \rightarrow_{\mathbb{P}} \beta^{-\lambda},
\end{equation}
\begin{equation}
\beta_{1} = \beta_{2} = \beta,
\end{equation}
\begin{equation}
\lim_{x \rightarrow \infty} e^{r x^{1/\lambda}} \mathbb{P}(X > x) =
\left \{
\begin{array}{ll}
0 & \forall~r < \beta ~\mbox{if}~ \beta > 0,\\
&\\
\infty & \forall~r > \beta ~\mbox{if}~ \beta < \infty.
\end{array}
\right.
\end{equation}
\end{theorem}

\vskip 0.2cm

{\it Proof}~~Since
\[
\max_{1 \leq k \leq n} \left(X_{k} \vee 1 \right) - \max_{1 \leq k \leq n} X_{k}
= \left \{
\begin{array}{ll}
0 & \mbox{if}~~X_{k} \geq 1 ~\mbox{for some}~1 \leq k \leq n,\\
&\\
1 - \max_{1 \leq k \leq n} X_{k}I_{\{X_{k} < 1 \}} & \mbox{if}~~X_{k} < 1 ~\mbox{for all}~1 \leq k \leq n,
\end{array}
\right.
\]
\[
0 \leq 1 - \max_{1 \leq k \leq n} X_{k}I_{\{X_{k} < 1 \}} \leq 1 - X_{1}I_{\{X_{1} < 1 \}}, ~\mbox{and}~
\lim_{n \rightarrow \infty} \frac{1 - X_{1}I_{\{X_{1} < 1 \}}}{(\log n)^{\lambda}} = 0 ~~\mbox{a.s.,}
\]
we have that
\[
\limsup_{n \rightarrow \infty} \frac{\max_{1 \leq k \leq n}X_{k}}{(\log n)^{\lambda}} =
\limsup_{n \rightarrow \infty} \frac{\max_{1 \leq k \leq n}\left(X_{k} \vee 1 \right)}{(\log n)^{\lambda}} ~~\mbox{a.s.}
\]
and
\[
\liminf_{n \rightarrow \infty} \frac{\max_{1 \leq k \leq n}X_{k}}{(\log n)^{\lambda}} =
\liminf_{n \rightarrow \infty} \frac{\max_{1 \leq k \leq n}\left(X_{k} \vee 1 \right)}{(\log n)^{\lambda}} ~~\mbox{a.s.}
\]
Note that
\[
\frac{\max_{1 \leq k \leq n}\left(X_{k} \vee 1 \right)}{(\log n)^{\lambda}} 
= \left(\frac{\max_{1 \leq k \leq n} \left(X_{k} \vee 1 \right)^{1/\lambda}}{\log n} \right)^{\lambda}
= \left(\frac{\log \max_{1 \leq k \leq n} e^{\left(X_{k} \vee 1 \right)^{1/\lambda}}}{\log n} \right)^{\lambda} ~~\forall~ n \geq 1
\]
and, for $y \geq e$
\[
\begin{array}{lll}
\mbox{$\displaystyle 
y^{r} \mathbb{P} \left(e^{(X \vee 1)^{1/\lambda}} > y \right)$}
& = & \mbox{$\displaystyle y^{r} \mathbb{P}\left( (X \vee 1)^{1/\lambda} > \log y \right)$} \\
&&\\
&=& \mbox{$\displaystyle y^{r} \mathbb{P}\left( X > (\log y)^{\lambda} \right)$}\\
&&\\
&=& \mbox{$\displaystyle e^{rx^{1/\lambda}} \mathbb{P} (X > x)$} ~~\mbox{(letting $\displaystyle y = e^{x^{1/\lambda}}$)}.  
\end{array}
\]
We thus see that
\[
\sup \left\{r \geq 0:~\lim_{x \rightarrow \infty} x^{r} \mathbb{P} \left(e^{(X \vee 1)^{1/\lambda}} > x \right) = 0 \right \}
= \sup \left\{r \geq 0:~\lim_{x \rightarrow \infty} e^{rx^{1/\lambda}} \mathbb{P} (X > x) = 0 \right \} = \beta_{1}
\]
and 
\[
\sup \left\{r \geq 0:~\liminf_{x \rightarrow \infty} x^{r} \mathbb{P} \left(e^{(X \vee 1)^{1/\lambda}} > x \right) = 0 \right \}
= \sup \left\{r \geq 0:~\liminf_{x \rightarrow \infty} e^{rx^{1/\lambda}} \mathbb{P} (X > x) = 0 \right \} = \beta_{2}
\]
and hence by Theorem 2.3
\[
\limsup_{n \rightarrow \infty} \frac{\max_{1 \leq k \leq n}X_{k}}{(\log n)^{\lambda}}
= \left(\limsup_{n \rightarrow \infty} 
\frac{\log \max_{1 \leq k \leq n} e^{\left(X_{k} \vee 1 \right)^{1/\lambda}}}{\log n} \right)^{\lambda}
= \beta_{1}^{-\lambda}~~\mbox{a.s.}
\]
and
\[
\liminf_{n \rightarrow \infty} \frac{\max_{1 \leq k \leq n}X_{k}}{(\log n)^{\lambda}}
= \left(\liminf_{n \rightarrow \infty} 
\frac{\log \max_{1 \leq k \leq n} e^{\left(X_{k} \vee 1 \right)^{1/\lambda}}}{\log n} \right)^{\lambda}
= \beta_{2}^{-\lambda}~~\mbox{a.s.}
\]
(i.e., (5.3) holds) and the statements (5.4), (5.5), (5.6), and (5.7) are equivalent. ~$\Box$

\vskip 0.3cm

\begin{remark}
For $0 < \lambda < \infty$ and $0 \leq \beta < \infty$, note that
\[
\mathbb{P} (X > x) 
= \exp \left(-\beta x^{1/\lambda} + \ln \left(e^{\beta x^{1/\lambda}} \mathbb{P}(X > x) \right) \right)~~\forall ~x > 0.
\]
We thus see that, if $0 < \lambda < \infty$ and $0 \leq \beta < \infty$, then (5.7) is equivalent to 
\[
\lim_{x \rightarrow \infty} \frac{\ln \left(e^{\beta x^{1/\lambda}} \mathbb{P}(X > x) \right)}{x^{1/\lambda}} = 0.
\]
Thus
\[
\lim_{n \rightarrow \infty} \frac{\max_{1 \leq k \leq n}X_{k}}{(\log n)^{\lambda}} 
= \beta^{-\lambda}~~\mbox{a.s. for some}~0 \leq \beta < \infty ~\mbox{and}~ 0 < \lambda < \infty
\]
if and only if there exists a function $H(\cdot)$ defined on $(0, \infty)$  such that
\[
\mathbb{P}(X > x) \sim \exp \left(- \beta x^{1/\lambda} + H(x) \right)~~\mbox{as}~~ x \rightarrow \infty
~~\mbox{and}~~\lim_{x \rightarrow \infty} \frac{H(x)}{x^{1/\lambda}} = 0.
\]
\end{remark}

\vskip 0.3cm

The following result is more general than that in Remark 5.1 provided $0 < \beta < \infty$.

\vskip 0.2cm

\begin{theorem}
Let $\{X_{n}; n \geq 1\}$ be a sequence of i.i.d. random variables drawn from 
the distribution function $F(\cdot)$ of the random variable $X$ such that
\begin{equation}
\mathbb{P}(X > x) \sim \exp\left(-\varphi(x) + H(x) \right) ~~\mbox{as}~~ x \rightarrow \infty
~~\mbox{and}~~\lim_{x \rightarrow \infty} \frac{H(x)}{\varphi(x)} = 0
\end{equation}
for some increasing and continuous function $\varphi(\cdot): (0, \infty) \rightarrow (0, \infty)$ and
some function $H(\cdot): (0, \infty)$ $\rightarrow (-\infty, \infty)$. If
\begin{equation}
\lim_{x \rightarrow \infty} \frac{\varphi^{-1}(x + {\it o}(x) )}{\varphi^{-1}(x)} = 1,
\end{equation}
then
\begin{equation}
\lim_{n \rightarrow \infty} \frac{\max_{1 \leq k \leq n}X_{k}}{\varphi^{-1}(\log n)} = 1 ~~\mbox{a.s.}
\end{equation}
\end{theorem}

\vskip 0.2cm

{\it Proof}~~Since $\varphi(\cdot):~(0, \infty) \rightarrow (0, \infty)$ is an increasing and continuous function,
\[
\begin{array}{lll}
\mbox{$\displaystyle 
e^{rx} \mathbb{P} \left(\varphi(X \vee 0) > x \right)$}
&=& \mbox{$\displaystyle
e^{rx} \mathbb{P} \left(X > \varphi^{-1}(x) \right)$}\\
&&\\
&=& \mbox{$\displaystyle 
e^{r\varphi(y)} \mathbb{P} (X > y)$ ~~(letting ~$\displaystyle y = \varphi^{-1}(x)$)~~$\displaystyle \forall ~x > 0$.}
\end{array}
\]
It thus follows from (5.8) that
\[
\lim_{x \rightarrow \infty} e^{rx} \mathbb{P} \left(\varphi(X \vee 0) > x \right)
= \left \{
\begin{array}{ll}
0 & \forall~r < 1,\\
&\\
\infty & \forall~r > 1
\end{array}
\right.
\]
and hence, by Theorem 5.1 with $\lambda = 1$ and $\beta = 1$, 
\[
\lim_{n \rightarrow \infty} \frac{\varphi\left(\max_{1 \leq k \leq n} \left(X_{k} \vee 0 \right) \right)}{\log n}
= \lim_{n \rightarrow \infty} \frac{\max_{1 \leq k \leq n} \varphi\left(X_{k} \vee 0 \right)}{\log n} = 1~~\mbox{a.s.};
\]
i.e., almost surely
\[
\varphi\left(\max_{1 \leq k \leq n} \left(X_{k} \vee 0 \right) \right) = \log n + {\it o} \left(\log n \right) ~~\mbox{as}~~
n \rightarrow \infty. 
\]
Thus (5.9) implies that almost surely
\[
\begin{array}{lll}
\mbox{$\displaystyle 
\max_{1 \leq k \leq n} \left(X_{k} \vee 0 \right)$}
&=& \mbox{$\displaystyle  \varphi^{-1}\left(\varphi\left(\max_{1 \leq k \leq n} \left(X_{k} \vee 0 \right) \right)\right)$}\\
&&\\
&=& \mbox{$\displaystyle
\varphi^{-1} \left(\log n + {\it o}(\log n) \right)$}\\
&&\\
&=& \mbox{$\displaystyle 
(1 + {\it o}(1)) \varphi^{-1}(\log n)~~\mbox{as}~~n \rightarrow \infty;$}
\end{array}
\]
i.e.,
\begin{equation}
\lim_{n \rightarrow \infty} \frac{\max_{1 \leq k \leq n} \left(X_{k} \vee 0 \right)}{\varphi^{-1}(\log n)} = 1~~\mbox{a.s.}
\end{equation}
Since
\[
\left| \max_{1 \leq k \leq n}X_{k} - \max_{1 \leq k \leq n} \left(X_{k} \vee 0 \right) \right |
\leq \left |\max_{1 \leq k \leq n} \left(X_{k} \wedge 0 \right) \right| \leq \left|\left(X_{1} \wedge 0 \right) \right|
~~\forall~n \geq 1,
\]
we see that (5.10) follows from (5.11).~~$\Box$

\vskip 0.3cm

Theorem 5.2 can be used to determine the asymptotic behavior very quickly for the maxima sequence
$\left\{\max_{1 \leq k \leq n} X_{k};~n \geq 1 \right \}$. This will be illustrated by the following 
three simple examples.

\vskip 0.3cm

\begin{example}
Let $\{X_{n}; n \geq 1\}$ be a sequence of i.i.d. random variables drawn from a standard normal random variable $X$.
It is well known that
\[
\mathbb{P}(X > x) = \frac{1}{\sqrt{2\pi}} \int_{x}^{\infty} e^{-z^{2}/2}dz \sim
\frac{1}{\sqrt{2\pi}x} e^{-x^{2}/2} = \exp\left(-\varphi(x) + H(x) \right)~~\mbox{as}~ x \rightarrow \infty
\]
where $\varphi(x) = x^{2}/2$ and  $H(x) = \ln \left(\sqrt{2\pi}x)^{-1}\right)$, $x > 0$.
Clearly, 
\[
\lim_{x \rightarrow \infty} \frac{H(x)}{\varphi(x)} = \lim_{x \rightarrow \infty}
\frac{\ln \left(\sqrt{2\pi}x)^{-1}\right)}{x^{2}/2} = 0,~~\varphi^{-1}(x) = \sqrt{2x} ~~\forall ~x > 0,
\]
and condition (5.9) is fulfilled. Thus, by Theorem 5.2, we have that 
\begin{equation}
\lim_{n \rightarrow \infty} \frac{\max_{1 \leq k \leq n}X_{k}}{\sqrt{2 \log n}} = 1~~\mbox{a.s.; i.e.,}~~
\lim_{n \rightarrow \infty} \frac{\max_{1 \leq k \leq n}X_{k}}{\sqrt{\log n}} = \sqrt{2}~~\mbox{a.s.}
\end{equation}
\end{example}

\vskip 0.3cm

\begin{remark}
Let $\{X_{n}; n \geq 1\}$ be a sequence of i.i.d. random variables drawn from a standard normal random variable $X$. 
We must point out that the stability properties for the maxima sequence $\left\{\max_{1 \leq k \leq n} X_{k};~n \geq 1 \right \}$
have been well studied by many authors. For example, Gnedenko (1943) proved that
\begin{equation}
\max_{1 \leq k \leq n} X_{k} - \sqrt{2 \log n} \rightarrow_{\mathbb{P}} 0
\end{equation}
which also yields (5.12) via Theorem 5.1. Goodman (1988) established a strong version of (5.13) in a Banach space setting. 
In particular, it follows from Theorem 2.1 of Goodman (1988) that
\[
\lim_{n \rightarrow \infty} 
\max_{1 \leq k \leq n} \left( \min_{-\sqrt{2 \log n} \leq x \leq \sqrt{2 \log n}} \left|X_{k} - x \right| \right) 
= 0~~\mbox{a.s.}
\]
and
\[
\lim_{n \rightarrow \infty} 
\max_{-\sqrt{2 \log n} \leq x \leq \sqrt{2 \log n}} \left(\min_{1 \leq k \leq n} \left|X_{k} - x \right| \right) = 0 ~~\mbox{a.s.}
\]
and hence
\[
\lim_{n \rightarrow \infty} \left(\max_{1 \leq k \leq n} X_{k} - \sqrt{2 \log n} \right) = 0 ~~\mbox{a.s.}
\]
\end{remark}

\vskip 0.3cm

\begin{example}
Let $\{X_{n}; n \geq 1\}$ be a sequence of i.i.d. random variables drawn from a Poisson random variable $X$ with parameter
$0 < \lambda < \infty$; i.e.,
\[
\mathbb{P}(X = m) = e^{-\lambda} \frac{\lambda^{m}}{m!},~~m = 0, 1, 2, ...
\]
Note that
\[
\lim_{m \rightarrow \infty} \frac{\mathbb{P}(X = m + 1)}{\mathbb{P}(X = m)} = 0
\]
and, by Stirling's formula,
\[
m! \sim \sqrt{2 \pi m} e^{-m} m^{m}~~\mbox{as}~~m \rightarrow \infty.
\]
Thus, as $x \rightarrow \infty$
\[
\begin{array}{lll}
\mbox{$\displaystyle 
\mathbb{P}(X > x)$} 
&\sim & \mbox{$\displaystyle
\mathbb{P}(X \geq [x] + 1)$}\\
&&\\
&\sim & \mbox{$\displaystyle 
e^{-\lambda} (2 \pi ([x]+1))^{-1/2} \lambda^{[x]+1} e^{-[x] -1} ([x] + 1)^{[x] + 1}$}\\
&&\\
&=& \mbox{$\displaystyle 
\exp \left(- ([x] + 1) \ln ([x] + 1) + ([x] + 1) (\ln \lambda - 1) - \frac{\lambda + \ln ([x]+1) + \ln (2\pi)}{2} \right)$} \\
&&\\
&=& \mbox{$\displaystyle
\exp\left(-\varphi(x) + H(x) \right)$,}
\end{array}
\]
where
\[
\varphi(x) = x \log x, ~~x > 0
\]
and
\[
H(x) = \left(x \log x - ([x] + 1) \ln ([x] + 1) \right) + ([x] + 1) (\ln \lambda - 1) - \frac{\lambda + \ln ([x]+1) + \ln (2\pi)}{2}, ~~x > 0.
\]
Note that 
\[
\varphi^{-1}(x) \sim \frac{x}{\log x}~~\mbox{and}~~\frac{x \log x - ([x] + 1) \ln ([x] + 1)}{x \log x} \rightarrow 0~~
\mbox{as}~~x \rightarrow \infty~~
\]
Thus it is easy to verify that all conditions of Theorem 5.2 are fulfilled with $\varphi(\cdot)$ and $H(\cdot)$.
Hence, by Theorem 4.2, we have that
\[
\lim_{n \rightarrow \infty} \frac{\max_{1 \leq k \leq n}X_{k}}{\frac{\log n}{\log \log n}} = 1~~\mbox{a.s.}
\]
\end{example}

\vskip 0.3cm

\begin{example}
Let $\{X_{n}; n \geq 1\}$ be a sequence of i.i.d. random variables drawn from a random variable $X$ with
\[
\mathbb{P}(X > x) = \exp \left(-e^{\zeta x^{\gamma}} + 1 \right), ~ x \geq 0
\]
where $\zeta > 0$ and $\gamma > 0$ are parameters. Let
\[
\varphi(x) = e^{\zeta x^{\gamma}} ~~\mbox{and}~~H(x) = 1, ~ x > 0.
\]
Then 
\[
\varphi^{-1}(x) = \left(\frac{\ln x}{\zeta} \right)^{1/\gamma}, ~~x > 1
\]
and all conditions of Theorem 5.2 are fulfilled with $\varphi(\cdot)$ and $H(\cdot)$. Thus, 
by Theorem 5.2, we have that
\[
\lim_{n \rightarrow \infty} \frac{\max_{1 \leq k \leq n} X_{k}}{\left(\frac{\log \log n}{\zeta} \right)^{1/\gamma}}
= 1 ~~\mbox{a.s.; i.e.,}~~\lim_{n \rightarrow \infty} \frac{\max_{1 \leq k \leq n} X_{k}}{\left(\log \log n \right)^{1/\gamma}}
= \zeta^{-1/\gamma} ~~\mbox{a.s.}
\]
\end{example}

\vskip 0.5cm

\noindent {\bf Acknowledgments}

\noindent The research of Deli Li was partially supported by a grant from the Natural Sciences 
and Engineering Research Council of Canada (Grant \#: RGPIN-2014-05428) and the research of 
Shuhua Zhang was partially supported by the National Natural Science Foundation of China (Grant 
\#: 91430108 and 11171251).

\newpage

{\bf References}

\begin{enumerate}

\item Barndorff-Nielsen, O.: On the limit behaviour of extreme order statistics. {\em Ann. Math. Statist.}
{\bf 34}, 992-1002 (1963)

\item Chandra, T.K.: {\em The Borel-Cantelli Lemma.} Springer, Heidelberg (2012)

\item Chow, Y. S., Teicher, H.: {\em Probability Theory: Independence, Interchangeability, Martingales., Third Ed.} 
Springer-Verlag, New York (1997)

\item Gnedenko, B.V.: Sur la distribution limite du terme maximum d'une s\'{e}rie al\'{e}atoire. {\em Ann. of Math.} (2)
{\bf 44}, 423-453 (1943)

\item Goodman, V.: Characteristics of normal samples. {\em Ann. Probab.} {\bf 16}, 1281-1290 (1988)

\item Tomkins, R.J.: Regular variation and the stability of maxima. {\em Ann. Probab.} {\bf 14}, 984-995 (1986)

\end{enumerate}

\end{document}